\documentclass{article}
\usepackage{amsfonts,amsmath,amssymb,amsthm,cite,enumerate,graphicx,latexsym,mathrsfs}
\usepackage{authblk}
\usepackage{hyperref}
\usepackage{float}
\makeatletter

\newcommand{\Rmnum}[1]{\expandafter\@slowromancap\romannumeral #1@}
\makeatother

\numberwithin{equation}{section}

\newtheorem{lemma}{Lemma}[section]

\newtheorem{remark}{Remark}[section]

\newtheorem{theorem}{Theorem}[section]

\newcommand\theref[1]{Theorem~\ref{#1}}
\newcommand\lemref[1]{Lemma~\ref{#1}}

\title{Global existence and stability of time-periodic solution to isentropic compressible Euler equations with source term}
\author{Huimin Yu\thanks{{e}-mail: hmyu@sdnu.edu.cn} \quad Xiaomin Zhang\thanks{{e}-mail: zxm15924687@163.com} \quad Jiawei Sun\thanks{Corresponding author {e}-mail: sunjiawei0122@163.com}
 \\ \small\textit{ Department of mathematics, Shandong Normal University, Jinan 250014 China}}

\begin{document}
\begin{sloppypar}
\date{}
\maketitle
\begin{center}
\begin{minipage}{130mm}{\small
\textbf{Abstract}:
In this paper, we study the initial-boundary value problem of one-dimensional isentropic compressible Euler equations with the source term $\beta\rho|u|^{\alpha}u$. By means of wave decomposition and the uniform a-priori estimates, we prove the global existence of smooth solutions under small perturbations around the supersonic Fanno flow. Then, by Gronwall's inequality, we get the smooth solution is time-periodic.
}
\end{minipage}
\end{center}
\section{Introduction}
\indent\indent  In this paper, we are concerned with isentropic compressible Euler equations with a nonlinear term:
\begin{equation}\label{a1}
\left\{\begin{aligned}
&\partial_{t}\rho+\partial_{x}(\rho u)=0,\\
&\partial_{t}(\rho u)+\partial_{x}(\rho u^{2}+p)=\beta\rho|u|^{\alpha}u,
\end{aligned}\right.
\quad(t,x)\in[0,+\infty)\times[0,L],
\end{equation}
where $\rho,u$ and $p$ are the density, velocity and pressure of gas, respectively. The pressure $p(\rho)$ is governed by $p(\rho)=a\rho^{\gamma}$, here the adiabatic exponent $\gamma>1$ and the parameter $a$ is scaled to unity for mathematical convenience. The sound speed $c\geq0$ is defined by $c^{2}=\partial p/\partial\rho$. And the term $\beta\rho|u|^{\alpha}u$ represents the friction with $\alpha,\beta\in\mathbb{R}$.\\
\indent In this paper, we assume the initial data are
\begin{align}\label{a2}
(\rho,u)^{\top}|_{t=0}=(\rho_{0}(x),u_{0}(x))^{\top}.
\end{align}
The boundary conditions are
\begin{align}\label{a3}
(\rho,u)^{\top}|_{x=0}=(\rho_{l}(t),u_{l}(t))^{\top}
\end{align}
and $\rho_{l}(t),u_{l}(t)$ are periodic functions with a period $P>0$, i.e. $$\rho_{l}(t+P)=\rho_{l}(t), u_{l}(t+P)=u_{l}(t). $$ In order to obtain the $C^{1}$ solution, the initial and boundary data should satisfy the following compatibility conditions at point $(0,0)$
\begin{equation}\label{a16}
\left\{
\begin{aligned}
&\rho_{l}'(0)+\rho_{0}'(0)u_{0}(0)+\rho_{0}(0)u_{0}'(0)=0,\\
&\rho_{l}'(0)u_{l}(0)+\rho_{l}(0)u_{l}'(0)+\rho_{0}'(0)u_{0}^{2}(0)
+2\rho_{0}(0)u_{0}(0)u_{0}'(0)\\
&\quad+p'_{0}(0)-\beta\rho_{0}(0)u_{0}^{\alpha+1}(0)=0,\\
&\rho_{0}(0)=\rho_{l}(0),u_{0}(0)
=u_{l}(0),
\end{aligned}\right.
\end{equation}
where
$$
p'_{0}(0)=\gamma\rho_{0}^{\gamma-1}(0)\rho'_{0}(0).
$$\\
\indent Because of the widespread application background, the compressible Euler equation with several kinds of source term have been studied extensively and there are fruitful results. For example, we can refer \cite{Hsiao1,Hsiao, Pan, Xin} for the research on the existence and stability of the small smooth solution, \cite{Chen, Ding, Huang,Hsiao2, Sui, Wang-Chen, Yin} for the singularity formation of smooth solution and the results on weak solution. In this paper, we are interested in the time-periodic solution of problem \eqref{a1}-\eqref{a3}. As far as we know, there are many works on the studies of time-periodic solutions of the partial differential equations such as the viscous fluids equations~\cite{Cai,Jin,Luo,Ma,Matsumura} and the hyperbolic conservation laws~\cite{G,Ohnawa,Takeno,Temple,Naoki}. All of the studies mentioned above discuss the time-periodic solutions which are derived by the time-periodic external forces or the piston motion. But there are few works on the time-periodic solutions of the hyperbolic conservations laws derived by the time-periodic boundary condition. In~\cite{Yuan}, Yuan studied time-periodic supersonic solutions for the isentropic compressible Euler equation (i.e. $\beta=0$) triggered by periodic supersonic boundary condition. For the quasilinear hyperbolic system with a more general time-periodic boundary conditions, Qu showed the existence and stability of the time-periodic solutions around a small neighborhood of $u\equiv0$ in~\cite{Qu}. Recently, Wei et al.~\cite{Wei} studied the global stability problem for supersonic flows in one dimensional compressible Euler equations with a friction term $-\mu\rho|u|u,\mu>0$. \\
\indent In this paper, we would like to show global existence and uniqueness of time-periodic supersonic solutions of initial-boundary value problem~\eqref{a1}-\eqref{a3} with the general friction term  $\beta\rho|u|^\alpha u$ by perturbing some supersonic Fanno flow. Different from \cite{Wei}, we consider \eqref{a1}-\eqref{a3} in the form of sound speed and fluid speed. Then the Fanno fluid are considered for some upstream with positive constants state $(c_-, u_-)$ at the left side. After analyzing the ODEs carefully, we get the maximal duct length $L_m$, exceed which the flow will get chock. Base on the supersonic Fanno flow, we prove the existence of time periodic solution by wave decomposition.\\
\indent The main results of this paper are:\\
\begin{theorem}\label{t3}
For any fixed non-sonic upstream state $(\rho_-, u_-)$ satisfying $0<u_-\neq \sqrt{\gamma}\rho_-^{{\gamma-1}\over 2}$, there exists a maximal duct length $L_m$, which only depend on $\alpha, \beta, \gamma$ and $(\rho_-, u_-)^{\top}$,  such that the steady solution $\tilde{V}=(\tilde{\rho}(x), \tilde{u}(x))^{\top}$ of problem \eqref{a1} exists in $[0, L_m]$ and keeps the upstream supersonic/subsonic state.
\end{theorem}

\begin{theorem}\label{t2}
Suppose the duct length $L<L_m$ and the upstream state $(\rho_-, u_-)$ is supersonic, i.e. $u_-> \sqrt{\gamma}\rho_-^{{\gamma-1}\over 2}$. Then there exists a $\varepsilon_{0}>0$ such that for any fixed $\varepsilon$ with $0<\varepsilon\leq\varepsilon_{0}$, if
\begin{align}
\|(\rho_0(x)-\tilde{\rho}(x),~u_0(x)-\tilde{u}(x))\|_{C^{1}([0,L])}<\varepsilon,\label{a19}\\
\|(\rho_l(t)-\rho_-, ~u_l(t)-u_-)\|_{C^{1}([0,+\infty))}<\varepsilon,\label{a20}
\end{align}
 then the mixed initial-boundary value problem~\eqref{a1} -\eqref{a3} have a unique $C^{1}$ solution $V=(\rho(t,x),u(t,x))^{\top}$ in the domain $E=\{(t,x)|t>0,x\in(0,L]\}$, satisfying
\begin{align*}
\|V-\tilde{V}\|_{C^{1}(E)}<C\varepsilon
\end{align*}
for some constant $C>0$ and
\begin{align*}
V(t+P,x)=V(t,x),\quad\forall t>T_{1},x\in[0,L],
\end{align*}
where $\tilde{V}=(\tilde{\rho}(x),\tilde{u}(x))^{\top}$ is the supersonic Fanno flow obtained in Theorem 1.1 and
\begin{align}\label{a21}
T_{1}=\max_{\substack{t\geq0,x\in[0,L]\\i=1,2}}\frac{L}{\lambda_{i}(V(t,x))}.
\end{align}
\end{theorem}
\begin{remark}
For the supersonic flow, the flow at $x=L$ is completely determined by the initial data at $x\in[0,L]$ and boundary conditions at $x=0$, so we only need to give the boundary condition at $x=0$.
\end{remark}

\indent The rest of the paper is organised as follows. In Section $2$, we construct the Fanno flow. In Section $3$, we present a reformulation of the problem by perturbing the solution around the supersonic Fanno flow and introduce a wave decomposition for the perturbed solution. In Section $4$, we prove global existence and uniqueness of solution under the help of uniform a-priori estimates. In Section $5$, we prove time-periodicity of solutions by the Gronwall's inequality.

%%%%%%%%%%%%%%%%%%%%%%%%%%%%%%%%%%%%%%%%%%%%%%%%%%%%%%%%%%%%%%%%%%%%%%%%%%%%%%%%%%%%%%%%%%%%%%%%%%%%%%%%%%%%%%%%%%%%%%%

\section{Fanno Flow}

\indent\indent Fanno flow refers to adiabatic flow through a constant area duct where the effect of friction $(i.e. \beta<0)$ is considered. The friction causes the flow properties to change along the duct. For the completeness of our results, we also consider the case $\beta>0$ in this section.\\
\indent We rewrite the initial-boundary problem~\eqref{a1}-~\eqref{a3} in terms of the sound speed $c=\sqrt{\gamma}\rho^{\frac{\gamma-1}{2}}$ and the fluid velocity $u$ as follows
\begin{equation}\label{a4}
\left\{\begin{aligned}
&c_{t}+c_{x}u+\frac{\gamma-1}{2}cu_{x}=0,\\
&u_{t}+uu_{x}+\frac{2}{\gamma-1}cc_{x}=\beta|u|^{\alpha}u,\\
&(c,u)^{\top}|_{t=0}=(c_{0}(x),u_{0}(x))^{\top},\\
&(c,u)^{\top}|_{x=0}=(c_{l}(t),u_{l}(t))^{\top},
\end{aligned}\right.
\end{equation}
where $c_{0}(x)=\sqrt{\gamma}\rho_{0}^{\frac{\gamma-1}{2}}(x),
c_{l}(t)=\sqrt{\gamma}\rho_{l}^{\frac{\gamma-1}{2}}(t)$.\\
\indent Now, we consider the positive solution $(\tilde{c},\tilde{u})^{\top}$ of the steady flow of system \eqref{a4} which satisfies
%use the above problem to construct steady state solution $(\tilde{c},\tilde{u})$, where $(\tilde{c},\tilde{u})$ is independent of $t$ and $\tilde{u}>0$. Then ~\eqref{a4} becomes
\begin{align}\label{a5}
\left\{
\begin{aligned}
&\tilde{c}'\tilde{u}+\frac{\gamma-1}{2}\tilde{c}\tilde{u}'=0,\\
&\tilde{u}\tilde{u}'+\frac{2}{\gamma-1}\tilde{c}\tilde{c}'
=\beta\tilde{u}^{1+\alpha},\\
&(\tilde{c},\tilde{u})^{\top}|_{x=0}=(c_{-},u_{-})^{\top},
\end{aligned}
\right.
\end{align}
%where $c_{-}=\frac{1}{P}\int_{0}^{P}c_{l}(t)dt=:\sqrt{\gamma}\rho_{-}^{\frac{\gamma-1}{2}},~u_{-}=\frac{1}{P}\int_{0}^{P}u_{l}(t)dt$.\\
where $u_{-}$ and $c_{-}$ are two positive constants.\\
\indent First, by $~\eqref{a5}_{1}$, we get
\begin{align}
\tilde{c}&=c_{-}u_{-}^{\frac{\gamma-1}{2}}\tilde{u}^{-\frac{\gamma-1}{2}}.\label{a6}
\end{align}
Substituting \eqref{a6} into $\eqref{a5}_{2}$, we have
\begin{align}
\tilde{u}^{-\alpha}\tilde{u}'-c_{-}^{2}u_{-}^{\gamma-1}\tilde{u}^{-\gamma-\alpha-1}\tilde{u}'
=\beta.\label{a7}
\end{align}
 We consider~\eqref{a7} by classifying $\alpha$ and $\beta$.\\\\
\textbf{Case 1:} $\alpha\neq1$ and $\alpha\neq-\gamma$.\\
\indent In this case, \eqref{a7} becomes
\begin{align}
\frac{1}{-\alpha+1}(\tilde{u}^{-\alpha+1})'+\frac{1}{\gamma+\alpha}c_{-}^{2}u_{-}^{\gamma-1}
(\tilde{u}^{-\gamma-\alpha})'=\beta.\label{a8}
\end{align}
Integrating~\eqref{a8} from $0$ to $x$, we get
\begin{align}\label{a9}
\frac{1}{-\alpha+1}\tilde{u}^{-\alpha+1}+\frac{1}{\gamma+\alpha}c_{-}^{2}u_{-}^{\gamma-1}
\tilde{u}^{-\gamma-\alpha}=\frac{1}{-\alpha+1}u_{-}^{-\alpha+1}+\frac{1}{\gamma+\alpha}
c_{-}^{2}u_{-}^{-1-\alpha}+\beta x.
\end{align}
Denote the left-hand-side function of \eqref{a9} as $h(s)$, i.e.
$$
h(s)=\frac{1}{-\alpha+1}s^{-\alpha+1}+\frac{1}{\gamma+\alpha}c_{-}^{2}u_{-}^{\gamma-1}
s^{-\gamma-\alpha},
$$
then we deduce
\begin{align*}
&h'(s)<0,\quad {\rm for~ }0<s<s_{c};\\
&h'(s)>0,\quad {\rm for~ }s>s_{c},
\end{align*}
where $s_{c}=c_{-}^{\frac{2}{\gamma+1}}u_{-}^{\frac{\gamma-1}{\gamma+1}}$.
This means that $h(s)$ gets its minimum at point $s=s_c$. On the other hand, from ~\eqref{a6}, we have $\tilde{c}=c_{-}^{\frac{2}{\gamma+1}}u_{-}^{\frac{\gamma-1}{\gamma+1}}$, when $\tilde{u}=s_{c}=c_{-}^{\frac{2}{\gamma+1}}u_{-}^{\frac{\gamma-1}{\gamma+1}}$. That is, the flow speed equals to the sound speed $(i.e. M=1)$ at the choked point $(s_c, h(s_c))$. See Figure 1 below.\\
\begin{figure}[H]
\centering
\includegraphics[width=9cm]{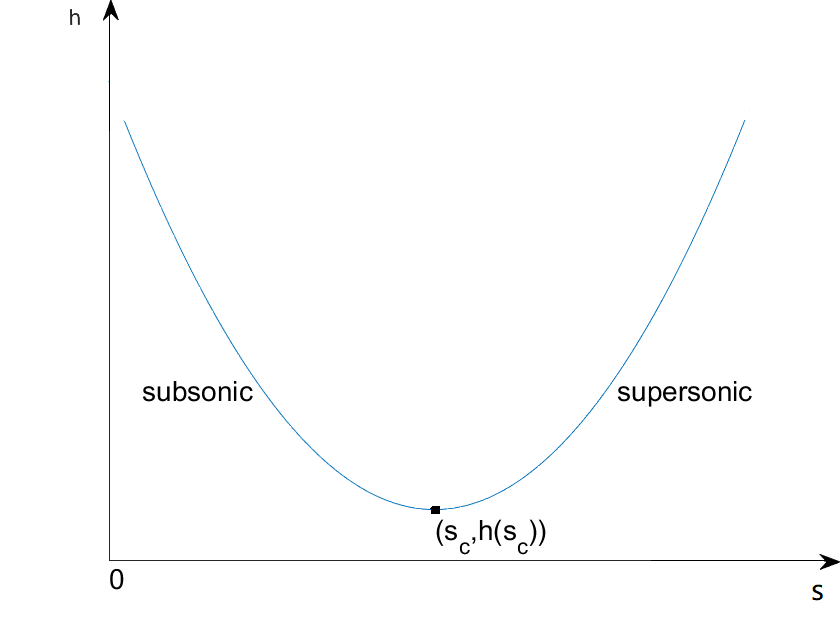}
\begin{center}
Figure 1
\end{center}
\end{figure}
\indent If $\beta>0$ and the upstream is supersonic (i.e. $u_->c_-$), $\tilde{u}$ is monotonically increasing by considering \eqref{a9} and $\tilde{u}> u_{-}$. By \eqref{a6}, $\tilde{c}$ is monotonically decreasing and $\tilde{c}< c_{-}$. Then, $\tilde{u}>\tilde{c}$. If $\beta>0$ and the upstream is subsonic (i.e. $u_-<c_-$), $\tilde{u}$ is monotonically decreasing and $\tilde{c}$ is monotonically increasing. Then $\tilde{u}<\tilde{c}$.\\
\indent When $\beta<0$, from \eqref{a9}, $h(s)$ decreases with respect to the length of the duct till arriving its minimum. Therefore, we can get the maximal length of the duct $L_{m}$ for a supersonic or subsonic flow before it gets choked, which is
\begin{align}\label{L1}
L_{m}=\frac{1}{\beta}\Big(\frac{1}{-\alpha+1}(s_{c}^{-\alpha+1}-u_{-}^{-\alpha+1})+\frac{1}{\gamma+\alpha}c_{-}^{2}(u_{-}^{\gamma-1}s_{c}^{-\gamma-\alpha}-u_{-}^{-1-\alpha})\Big).
\end{align}
\textbf{Case 2:} $\alpha=1$ or $\alpha=-\gamma$.\\
\indent Now, \eqref{a7} is turned into
\begin{align}
(\ln\tilde{u})'+\frac{1}{\gamma+1}c_{-}^{2}u_{-}^{\gamma-1}
(\tilde{u}^{-\gamma-1})'=\beta, \quad {\rm for} ~~\alpha=1,\label{a10}
\end{align}
and
\begin{align}
\frac{1}{\gamma+1}(\tilde{u}^{\gamma+1})'-c_{-}^{2}u_{-}^{\gamma-1}
(\ln\tilde{u})'=\beta, \quad {\rm for} ~~\alpha=-\gamma. \label{a12}
\end{align}
Integrating ~\eqref{a10} and~\eqref{a12} from $0$ to $x$, we get
\begin{align}
\ln\tilde{u}+\frac{1}{\gamma+1}c_{-}^{2}u_{-}^{\gamma-1}
\tilde{u}^{-\gamma-1}=\ln u_{-}+\frac{1}{\gamma+1}c_{-}^{2}u_{-}^{-2}+\beta x, \quad {\rm for} ~~\alpha=1, \label{a11}
\end{align}
and
\begin{align}
\frac{1}{\gamma+1}\tilde{u}^{\gamma+1}-c_{-}^{2}u_{-}^{\gamma-1}\ln\tilde{u}
=\frac{1}{\gamma+1}u_{-}^{\gamma+1}-c_{-}^{2}u_{-}^{\gamma-1}\ln u_{-}+\beta x, \quad {\rm for}~~\alpha=-\gamma. \label{a13}
\end{align}
Define
$$
f(s)=\ln s+\frac{1}{\gamma+1}c_{-}^{2}u_{-}^{\gamma-1}
s^{-\gamma-1},
$$
and
$$
g(s)=\frac{1}{\gamma+1}s^{\gamma+1}-c_{-}^{2}u_{-}^{\gamma-1}\ln s.
$$
The functions $f(s)$ and $g(s)$ get their minimums at point $s=s_c=c_{-}^{\frac{2}{\gamma+1}}u_{-}^{\frac{\gamma-1}{\gamma+1}}$.
Furthermore, we get the maximal length of the duct $L_{m}$ :
\begin{align}\label{L2}
L_{m}=\frac{1}{\beta}\Big(\frac{1}{\gamma+1}c_{-}^{2}(u_-^{\gamma-1}s_c^{-\gamma-1}-u_-^{-2})+
\ln\frac{s_c}{u_{-}}\Big), {~~\rm for~~} \alpha=1
\end{align}
and
\begin{align}\label{L3}
L_{m}=\frac{1}{\beta}\Big(\frac{1}{\gamma+1}(s_c^{\gamma+1}-u_{-}^{\gamma+1})-c_{-}^{2}
u_{-}^{\gamma-1}\ln\frac{s_c}{u_{-}}\Big), {~~\rm for~~} \alpha=-\gamma.
\end{align}
\indent We can get the similar results with the case 1, we omit the details here.\\
\indent From the above discussion, we have the following Lemma.
\begin{lemma}\label{t1}
If $u_->0, c_->0$ and the duct length $L< L_m$, where $L_m$ is a positive constant only depending on $\alpha, \beta, \gamma, c_-$ and $u_-$ (See~\eqref{L1},\eqref{L2},\eqref{L3}), then the Cauchy problem~\eqref{a5} admits a unique smooth positive solution $(\tilde{c}(x),\tilde{u}(x))^{\top}$ which satisfies the following properties:
\begin{enumerate}[1)]
\item  $0<\tilde{u}(x)<u_-<c_-<\tilde{c}(x)$,~~~ if  $\beta>0$ and $c_->u_-$;
\item  $0<\tilde{c}(x)<c_-<u_-<\tilde{u}(x)$,~~~ if  $\beta>0$ and $c_-<u_-$;
\item  $0<u_-<\tilde{u}(x)<\tilde{c}(x)<c_-$,~~~ if  $\beta<0$ and $c_->u_-$;
\item  $0<c_-<\tilde{c}(x)<\tilde{u}(x)<u_-$,~~~ if  $\beta<0$ and $c_-<u_-$.
\end{enumerate}
\end{lemma}
\indent This result means that a subsonic flow entering a duct with friction $(\beta<0)$ will have an increase in its Mach number until the flow is choked at $M=1$, i.e. $\tilde{u}=\tilde{c}$. Conversely, the Mach number of a supersonic flow will decrease until the flow is choked. However, if a flow entering a duct with acceleration $(\beta>0)$, the Mach number of a subsonic flow will decrease and the Mach number of a supersonic flow will increase (i.e. accelerating the initial subsonic or supersonic state). It is worthy to be pointed out that the theoretic calculations are consistent with its experiment. Different from the calculations in \cite{Wei},  where the authors consider a differential equation that relates the change in Mach number with respect to the length of the duct ${dM}\over {dx}$, we rewrite the dominating equations by the relations between the sound speed and flow speeds. Fortunately, the resulting equations can be decoupled easily. Therefore, we can show the maximal duct length which makes the flow choke assuming the upstream Mach number is supersonic or subsonic. Thus by~\lemref{t1} and $\tilde{c}=\sqrt{\gamma}\tilde{\rho}^{\frac{\gamma-1}{2}}$, we can directly get~\theref{t3}.\\
 \indent From the result, we observe that no matter what the constant number $\alpha$ is, the Mach number of a supersonic (subsonic) flow will increase (decrease) when $\beta>0$. While when $\beta<0$, the Mach number of a supersonic flow will decrease until the flow is choked; conversely, a subsonic flow will have an increase in its Mach number until the flow is choked.\\

\section{Reformulation of Problem and Wave Decomposition}\label{s2}

\indent\indent For the supersonic flow, we should have $u>0$. Then, we can write the system~\eqref{a1} as
\begin{align}\label{b1}
\left\{
\begin{aligned}
&\rho_{t}+\rho_{x}u+\rho u_{x}=0,\\
&u_{t}+uu_{x}+\gamma\rho^{\gamma-2}\rho_{x}=\beta u^{\alpha+1}.
\end{aligned}\right.
\end{align}
\indent Letting
\begin{align}\label{b2}
\rho(t,x)=\bar{\rho}(t,x)+\tilde{\rho}(x),\quad u(t,x)=\bar{u}(t,x)
+\tilde{u}(x),
\end{align}
where $(\bar{\rho}(t,x),\bar{u}(t,x))^{\top}$ is the perturbation of the supersonic Fanno flow. Substituting ~\eqref{b2} into ~\eqref{b1}, we get
\begin{align}\label{b3}
\left\{
\begin{aligned}
&\bar{\rho}_{t}+\bar{\rho}_{x}u+\rho \bar{u}_{x}+\tilde{\rho}'\bar{u}+\bar{\rho}\tilde{u}'
+\tilde{\rho}'\tilde{u}+\tilde{\rho}\tilde{u}'=0,\\
&\bar{u}_{t}+u\bar{u}_{x}+\bar{u}\tilde{u}'
+\tilde{u}\tilde{u}'+\gamma\rho^{\gamma-2}\bar{\rho}_{x}
+\gamma\rho^{\gamma-2}\tilde{\rho}'=\beta(\bar{u}+\tilde{u})^{\alpha+1}.
\end{aligned}\right.
\end{align}
 Moreover, the system ~\eqref{b3} can be further written into
\begin{align}\label{b4}
\left\{
\begin{aligned}
&\bar{\rho}_{t}+\bar{\rho}_{x}u+\rho \bar{u}_{x}=-\tilde{\rho}'\bar{u}-\bar{\rho}\tilde{u}',\\
&\bar{u}_{t}+u\bar{u}_{x}+\gamma\rho^{\gamma-2}\bar{\rho}_{x}
=-F(\rho,\tilde{\rho})\bar{\rho}\tilde{\rho}'-\bar{u}\tilde{u}'
-G(u,\tilde{u})\bar{u},
\end{aligned}\right.
\end{align}
where $F(\rho,\tilde{\rho})\bar{\rho}=\gamma(\rho^{\gamma-2}-\tilde{\rho}^{\gamma-2})$, $G(u,\tilde{u})\bar{u}=-\beta[u^{\alpha+1}-\tilde{u}^{\alpha+1}]$ and $F(\rho,\tilde{\rho})$ and $G(u,\tilde{u})$ can be taken the following expressions
\begin{align}
F(\rho,\tilde{\rho})=\gamma(\gamma-2)\int_{0}^{1}(\theta\bar{\rho}+\tilde{\rho})^{\gamma-3}d\theta,\quad
G(u,\tilde{u})=-\beta(\alpha+1)\int_{0}^{1}(\theta\bar{u}+\tilde{u})^{\alpha}d\theta. \label{Exps1}
\end{align}
\indent We also consider the perturbations of the initial and boundary conditions. The initial data is reformulated as
\begin{align}\label{a14}
t=0:\left\{
\begin{aligned}
\rho_{0}(x)=\bar{\rho}_{0}(x)+\tilde{\rho}(x),\quad x\in[0,L],\\\
u_{0}(x)=\bar{u}_{0}(x)+\tilde{u}(x),\quad x\in[0,L],
\end{aligned}\right.
\end{align}
where $L<L_{m}$,
and boundary condition is
\begin{align}\label{a15}
x=0:\left\{
\begin{aligned}
&\rho_{l}(t)=\bar{\rho}_{l}(t)+\tilde{\rho}(0),\quad t\geq0,\\
&u_{l}(t)=\bar{u}_{l}(t)+\tilde{u}(0),\quad t\geq0,
\end{aligned}\right.
\end{align}
where $\bar{\rho}_{0},\bar{u}_{0}, \bar{\rho}_{l},\bar{u}_{l}$ are $C^{1}$ functions.\\

\indent Let $\bar{V}=(\bar{\rho},\bar{u})^{\top}$, the system~\eqref{b4} can be rewritten as the following quasi-linear form
\begin{equation}\label{b5}
\bar{V}_{t}+A(V)\bar{V}_{x}+D(\tilde{V})\bar{V}=0
\end{equation}
with the initial data
\begin{align}
\bar{V}|_{t=0}=\bar{V}_{0}=(\bar{\rho}_{0},\bar{u}_{0})^{\top},\label{1in1}
\end{align}
and the boundary condition
\begin{align}
V|_{x=0}=V_{l}&=(\rho_{l},u_{l})^{\top},\label{1in2}
\end{align}
where $V(t,x)=\bar{V}(t,x)+\tilde{V}(x)$, and
$$
A(V)=\left(
 \begin{matrix}
 u & \rho\\
 \gamma\rho^{\gamma-2} & u
 \end{matrix}
\right),\quad
D(\tilde{V})=\left(
 \begin{matrix}
 \tilde{u}' & \tilde{\rho}'\\
 F(\rho,\tilde{\rho})\tilde{\rho}' & \tilde{u}'+G(u,\tilde{u})
 \end{matrix}
\right).
$$
\indent We next introduce a wave decomposition of the solution $\bar{V}$ to the system~\eqref{b5}. We can easily get the following two eigenvalues of the coefficient matrix $A(V)$
$$
\lambda_{1}(V)=u-c,\quad\lambda_{2}(V)=u+c,
$$
where $c=\sqrt{\gamma}\rho^{\frac{\gamma-1}{2}}$. The corresponding two right eigenvectors $r_{i}, i=1,2$ are
\begin{align}\label{b6}
r_{1}(V)=\frac{1}{\sqrt{\rho^{2}+c^{2}}}(\rho,-c)^{\top},\quad
r_{2}(V)=\frac{1}{\sqrt{\rho^{2}+c^{2}}}(\rho,c)^{\top}.
\end{align}
The left eigenvectors $l_{i}(V), i=1,2$ are determined by
\begin{align}\label{b7}
l_{i}(V)r_{j}(V)\equiv\delta_{ij},\quad r_{i}^{\top}(V)r_{i}(V)\equiv1,\quad i,j=1,2,
\end{align}
where $\delta_{ij}$ stands for the Kronecker's symbol. Then, $l_{i}, i=1,2$ have the following expressions
\begin{align}
l_{1}(V)=\frac{\sqrt{\rho^{2}+c^{2}}}{2}(\rho^{-1},-c^{-1}),\quad l_{2}(V)=\frac{\sqrt{\rho^{2}+c^{2}}}{2}(\rho^{-1},c^{-1}), \label{Leg1}
\end{align}
 which have the same regularity as $r_{i}(V)$.\\
\indent Let
\begin{align}\label{b8}
m_{i}=l_{i}(V)\bar{V},\quad n_{i}=l_{i}(V)\bar{V}_{x},\quad m=(m_{1},m_{2})^{\top}
,\quad n=(n_{1},n_{2})^{\top},
\end{align}
then
\begin{align}
\bar{V}&=\sum_{k=1}^{2}m_{k}r_{k}(V),\quad
\frac{\partial\bar{V}}{\partial x}=\sum_{k=1}^{2}n_{k}r_{k}(V),\label{Sdec1}\\
\frac{\partial\bar{V}}{\partial t}&=-D(\tilde{V})\bar{V}-\sum_{k=1}^{2}\lambda_{k}(V)n_{k}r_{k}(V). \label{Sdec2}
\end{align}
Thus, we have
\begin{equation}\label{b9}
\begin{split}
\frac{d\bar{V}}{d_{i}t}&=\frac{\partial\bar{V}}{\partial t}+\lambda_{i}(V)\frac{\partial\bar{V}}{\partial x}\\
&=\sum_{k=1}^{2}(\lambda_{i}(V)-\lambda_{k}(V))n_{k}r_{k}(V)
-D(\tilde{V})\bar{V}.
\end{split}
\end{equation}
\indent By~\eqref{b7}-~\eqref{b9}, one has
\begin{equation}\label{b10}
\begin{split}
\frac{d m_{i}}{d_{i}t}=&\frac{\partial m_{i}}{\partial t}+\lambda_{i}(V)\frac{\partial m_{i}}{\partial x}\\
=&\sum_{j,k=1}^{2}\Psi_{ijk}(V)n_{j}m_{k}+\sum_{j,k=1}^{2}\tilde{\Psi}_{ijk}(V)
m_{j}m_{k}-\sum_{k=1}^{2}\tilde{\tilde{\Psi}}_{ik}(V)m_{k},
\end{split}
\end{equation}
where
\begin{align}
\Psi_{ijk}(V)&=(\lambda_{j}(V)-\lambda_{i}(V))l_{i}(V)r_{j}(V)\cdot\nabla_{V}r_{k}(V), \label{mep1}\\
\tilde{\Psi}_{ijk}(V)&=l_{i}(V)D(\tilde{V})r_{j}(V)\cdot\nabla_{V}r_{k}(V), \label{mep2}\\
\tilde{\tilde{\Psi}}_{ik}(V)&=\lambda_{i}(V)
l_{i}(V)\tilde{V}'\cdot\nabla_{V}r_{k}(V)+l_{i}(V)D(\tilde{V}) r_{k}(V). \label{mep3}
\end{align}
\indent Similarly, using~\eqref{b5} and ~\eqref{b7}-~\eqref{b9}, we also get
\begin{equation}\label{b11}
\begin{split}
\frac{d n_{i}}{d_{i}t}=&\frac{\partial n_{i}}{\partial t}+\lambda_{i}(V)\frac{\partial n_{i}}{\partial x}\\
=&\sum_{j,k=1}^{2}\Phi_{ijk}(V)n_{j}n_{k}+\sum_{j,k=1}^{2}\tilde{\Phi}_{ijk}(V)n_{k}-\sum_{k=1}^{2}l_{i}(V)D_{x}(\tilde{V})r_{k}(V)m_{k},
\end{split}
\end{equation}
where the term $D_{x}(\tilde{V})$ makes sense if $\tilde{V}$ is a $C^{2}$ function, and
\begin{align}
\Phi_{ijk}(V)=&(\lambda_{j}(V)-\lambda_{k}(V))l_{i}(V)r_{j}(V)\cdot\nabla_{V}r_{k}(V)\notag\\
&-r_{j}(V)\cdot\nabla_{V}\lambda_{k}(V)\delta_{ik},\label{nep1}\\
\tilde{\Phi}_{ijk}(V)=&-\lambda_{k}(V)l_{i}(V)\tilde{V}'\cdot\nabla_{V}r_{k}(V)
+l_{i}(V)D(\tilde{V})r_{j}(V)\cdot\nabla_{V}r_{k}(V)m_{j}(V)\notag\\
&-l_{i}(V)D(\tilde{V})r_{k}(V)-\tilde{V}'\cdot\nabla_{V}\lambda_{k}(V)\delta_{ik}. \label{nep2}
\end{align}
\indent For later use, we rewrite the system ~\eqref{b4} by exchanging the variable $t$ and $x$ as follows
\begin{equation*}
\bar{V}_{x}+A^{-1}(V)\bar{V}_{t}+A^{-1}(V)D(\tilde{V})\bar{V}=0.
\end{equation*}
Denote $\hat{\lambda}_{i}(V), i=1,2$ are eigenvalues of the matrix $A^{-1}(V)$, $\hat{l}_{i}(V), i=1,2$ and $\hat{r}_{i}(V), i=1,2$ are the left and right eigenvectors respectively. They can be determined in terms of $\lambda_{i}(V), r_{i}(V)$ and $l_{i}(V)$ as follows
\begin{align}
\hat{\lambda}_{i}(V)=\lambda_{i}(V)^{-1},\quad \hat{r}_{i}(V)=r_{i}(V),\quad \hat{l}_{i}(V)=l_{i}(V). \label{Evf1}
\end{align}
Therefore, $\hat{r}_{i}(V)$ and $\hat{l}_{i}(V)$ also satisfy~\eqref{b7}.

\indent Let
\begin{align}
\hat{m}_{i}=\hat{l}_{i}(V)\bar{V},\quad\hat{n}_{i}=\hat{l}_{i}(V)\bar{V}_{t},\quad
\hat{m}=(\hat{m}_{1},\hat{m}_{2})^{\top},\quad\hat{n}=(\hat{n}_{1},\hat{n}_{2})^{\top}. \label{hsdec1}
\end{align}
By applying the similar arguments as in~\eqref{b10}-~\eqref{nep2}, we can get
\begin{equation}\label{b12}
\begin{split}
\frac{d \hat{m}_{i}}{d_{i}x}=&\frac{\partial\hat{m}_{i}}{\partial x}+\hat{\lambda}_{i}(V)\frac{\partial\hat{m}_{i}}{\partial t}\\
=&\sum_{j,k=1}^{2}\hat{\Psi}_{ijk}(V)\hat{n}_{j}\hat{m}_{k}
+\sum_{j,k=1}^{2}\hat{\tilde{\Psi}}_{ijk}(V)\hat{m}_{j}\hat{m}_{k}
-\sum_{k=1}^{2}\hat{\tilde{\tilde{\Psi}}}_{ik}(V)\hat{m}_{k}
\end{split}
\end{equation}
with
\begin{align}
\hat{\Psi}_{ijk}(V)=&(\hat{\lambda}_{j}(V)-\hat{\lambda}_{i}(V))
\hat{l}_{i}(V)\hat{r}_{j}(V)\cdot\nabla_{V}\hat{r}_{k}(V),\label{hmep1}\\
\hat{\tilde{\Psi}}_{ijk}(V)=&\hat{\lambda}_{i}(V)\hat{l}_{i}(V)
D(\tilde{V})\hat{r}_{j}(V)\cdot\nabla_{V}\hat{r}_{k}(V),\label{hmep2}\\
\hat{\tilde{\tilde{\Psi}}}_{ik}(V)=&\hat{l}_{i}(V)\tilde{V}'\cdot\nabla_{V}
\hat{r}_{k}(V)+\hat{\lambda}_{i}(V)\hat{l}_{i}(V)D(\tilde{V})\hat{r}_{k}(V),\label{hmep3}
\end{align}
and
\begin{equation}\label{b13}
\begin{split}
\frac{d \hat{n}_{i}}{d_{i}x}=&\frac{\partial\hat{n}_{i}}{\partial x}+\hat{\lambda}_{i}(V)\frac{\partial\hat{n}_{i}}{\partial t}\\
=&\sum_{j,k=1}^{2}\hat{\Phi}_{ijk}(V)\hat{n}_{j}\hat{n}_{k}
+\sum_{j,k=1}^{2}\hat{\tilde{\Phi}}_{ijk}(V)\hat{n}_{k}-\sum_{k=1}^{2}\hat{l}_{i}(V)(A^{-1}(V)D(\tilde{V}))_{t}\hat{r}_{k}(V)\hat{m}_{k}(V)
\end{split}
\end{equation}
with
\begin{align*}
\hat{\Phi}_{ijk}(V)=&(\hat{\lambda}_{j}(V)-\hat{\lambda}_{k}(V))\hat{l}_{i}(V)\hat{r}_{j}(V)\cdot\nabla_{V}
\hat{r}_{k}(V)-\hat{r}_{j}(V)\cdot\nabla_{V}\hat{\lambda}_{k}(V)\delta_{ik},\\
\hat{\tilde{\Phi}}_{ijk}(V)=&-\hat{l}_{i}(V)\tilde{V}'\cdot\nabla_{V}\hat{r}_{k}(V)+
\hat{\lambda}_{i}(V)\hat{l}_{i}(V)D(\tilde{V})\hat{r}_{j}(V)\cdot\nabla_{V}\hat{r}_{k}(V)\hat{m}_{j}(V)\\
&-\hat{\lambda}_{i}(V)\hat{l}_{i}(V)D(\tilde{V})\hat{r}_{k}(V).
\end{align*}
\indent We also provide the wave decomposition of the initial and boundary data as follows
\begin{align}\label{ib8}
m_{0}=(m_{10},m_{20})^{\top}
,\quad n_{0}=(n_{10},n_{20})^{\top}
\end{align}
with
$$
m_{i0}=l_{i}(V_{0})\bar{V}_{0},\quad n_{i0}=l_{i}(V_{0})\bar{V}'_{0},
$$
and
\begin{align}
\hat{m}_{l}=(\hat{m}_{1l},\hat{m}_{2l})^{\top},\quad\hat{n}_{l}=(\hat{n}_{1l},\hat{n}_{2l})^{\top}. \label{ihsdec1}
\end{align}
with
$$
\hat{m}_{il}=\hat{l}_{i}(V_{l})\bar{V}_{l},\quad\hat{n}_{il}=\hat{l}_{i}(V_{l})\bar{V}'_{l},
$$
where $\bar{V}_{0}$ and $\bar{V}_{l}$ are defined by~\eqref{1in1} and~\eqref{1in2} respectively, and
\begin{align}
V_{0}&=(\rho_{0},u_{0})^{\top},\quad \bar{V}'_{0}=(\bar{\rho}'_{0},\bar{u}'_{0})^{\top},\label{in1}\\
V_{l}&=(\rho_{l},u_{l})^{\top},\quad \bar{V}'_{l}=(\bar{\rho}'_{l},\bar{u}'_{l})^{\top}.\label{in2}
\end{align}
\section{Existence of Global Solutions}\label{s3}
\indent\indent In this section, we will prove the existence of global solution $\bar{V}=(\bar{\rho}(t,x),\bar{u}(t,x))^{\top}$ to the initial-boundary value problem~\eqref{b5} and~\eqref{1in1}-\eqref{1in2} in the domain $E=\{(t,x)|t>0,x\in(0,L]\}$.\\
\indent The local existence and uniqueness of the $C^{1}$ solution to the mixed initial-boundary value problem~\eqref{b5} and~\eqref{1in1}-\eqref{1in2} is guaranteed by the classical theory in~\cite{Yu}, which can be extended globally in terms of a uniform a-priori estimate of the global $C^{1}$ solutions (see~\cite{Wei, Li, Zhou, Wang,Rao,Jin}).

\indent Next we will establish a uniform a-priori estimate of the classical solution to help us to extend globally the local solution. Let us first give the following assumption
\begin{align}
|m_{i}(t,x)|,\, |n_{i}(t,x)|\leq C\varepsilon,\quad\forall i=1,2, \quad (t,x)\in E  \label{c1}
\end{align}
for a suitably small positive constant $\varepsilon$, which will be determined later.

From\eqref{b6},~\eqref{Sdec1}and~\eqref{c1}, we have
\begin{align}
|\bar{V}(t,x)|,\,|\frac{\partial\bar{V}}{\partial x}(t,x)|\leq C\varepsilon,\quad\forall (t,x)\in E. \label{Casp2}
\end{align}
Combining~\theref{t1} with~\eqref{Casp2}, we obtain the following results. The details of the proof are omitted here.
\begin{lemma}\label{Casp3}
For sufficiently small $\varepsilon$, it holds that
\begin{align}
&|D(\tilde{V})(t,x)|, |\partial_{x}D(\tilde{V})(t,x)|, |\nabla_{V}r_{i}(V)(t,x)|, |\tilde{V}'|, T_{1}\leq C,\label{Casp4}\\
&C^{-1}\leq |\lambda_{i}(V)(t,x)|, |\nabla_{V}\hat{\lambda}_{i}(V)(t,x)|, |l_{i}(V)(t,x)|\leq C, \label{Casp5}\\
&|\frac{\partial\bar{V}}{\partial t}(t,x)|, |\partial_{t}A^{-1}(V)(t,x)|, |\partial_{t}D(\tilde{V})(t,x)|\leq C\varepsilon \label{Casp6}
\end{align}
for any $(t,x)\in E$, where the positive constant $C$ only depends on $c_{-}, u_{-}, \tilde{c}(L)$, $\tilde{u}(L), \gamma, \alpha$ and $\beta$.
\end{lemma}
We observe from~\eqref{Casp2} and~\eqref{Casp5} that it suffices to prove~\eqref{c1} for a uniform a-priori estimate of the global $C^{1}$ solution.

\indent Write $x=x_{i}^{*}(t), i=1,2$ as the characteristic curve of $\lambda_{i}$ passing a point $(0,0)$, which satisfy
\begin{align*}
\frac{d x_{i}^{*}(t)}{dt}=\lambda_{i}(V(t,x_{i}^{*}(t))),\quad x_{i}^{*}(0)=0.
\end{align*}
Noting that $x=x_{2}^{*}(t)$ lies below $x=x_{1}^{*}(t)$ since $\lambda_{2}(V)>\lambda_{1}(V)$.\\
\indent We divide the region $E$ into three small regions and discuss the uniform a-priori estimate of the classical solutions in each small region separately.\\
\textbf{Region 1:} the region $E_{1}=\{(t,x)|0\leq t\leq T_{1}, 0\leq x\leq L, x\geq x_{2}^{*}(t)\}$.\\
\indent For any point $(t,x)\in E_{1}$, integrating the i-th equation in ~\eqref{b10} along the i-characteristic curve with respect to $\tau$ from $0$ to $t$ which intersects the $x$-axis at a point $(0,b_{i})$, we obtain from~\eqref{b10}, \eqref{mep1}-\eqref{mep3}, \eqref{c1}, \eqref{Casp4} and~\eqref{Casp5} that
\begin{equation}\label{c2}
\begin{aligned}
|m_{i}(t,x(t))|\leq&|m_{i}(0,b_{i})|+\int_{0}^{t}\sum_{j,k=1}^{2}|\Psi_{ijk}(V)n_{j}m_{k}|d\tau\\
&+\int_{0}^{t}\sum_{j,k=1}^{2}|\tilde{\Psi}_{ijk}(V)m_{j}m_{k}|d\tau+\int_{0}^{t}\sum_{k=1}^{2}|\tilde{\tilde{\Psi}}_{ik}(V)m_{k}|d\tau\\
\leq&|m_{i0}(b_{i})|+C\int_{0}^{t}|m(\tau,x(\tau))|d\tau.
\end{aligned}
\end{equation}
\indent Applying the same procedures as above for~\eqref{b11}, from~\eqref{nep1}, \eqref{nep2}, \eqref{c1}, \eqref{Casp4} and~\eqref{Casp5}, we have
\begin{equation}\label{c3}
\begin{split}
|n_{i}(t,x(t))|\leq&|n_{i}(0,b_{i})|+\int_{0}^{t}\sum_{j,k=1}^{2}|\Phi_{ijk}(V)
n_{j}n_{k}|d\tau\\
&+\int_{0}^{t}\sum_{j,k=1}^{2}|\tilde{\Phi}_{ijk}(V)n_{k}|d\tau+\int_{0}^{t}\sum_{k=1}^{2}|l_{i}(V)D_{x}(\tilde{V})r_{k}(V)m_{k}|d\tau\\
\leq&|n_{i0}(b_{i})|+C(\int_{0}^{t}|n(\tau,x(\tau))|d\tau+\int_{0}^{t}|m(\tau,x(\tau))|d\tau).
\end{split}
\end{equation}
\indent Putting~\eqref{c2}-\eqref{c3} together, summing up $i=1,2$ and applying the Gronwall's inequality, we have
\begin{equation}\label{c4}
|m(t,x)|+|n(t,x)|\leq (\|m_{0}\|_{C^{0}([0,L])}+\|n_{0}\|_{C^{0}([0,L])})(1+CT_{1}).
\end{equation}
Because of the arbitrariness of $(t,x)\in E_{1}$ and the boundedness of $T_{1}$ in~\eqref{Casp4}, we obtain from~\eqref{c4} that
\begin{align}
\max_{(t,x)\in E_{1}}|m(t,x)|+|n(t,x)|\leq C(\|m_{0}\|_{C^{0}([0,L])}+\|n_{0}\|_{C^{0}([0,L])}). \label{Apri1}
\end{align}
\textbf{Region 2:} the region $E_{2}=\{(t,x)|t\geq0, 0\leq x\leq L, 0\leq x\leq x_{1}^{*}(t)\}$.\\
\indent For any point $(t,x)\in E_{2}$, integrating in~\eqref{b12} with respect to $x$ along the $i$-th characteristic curve, which is assumed to intersect the $t$-axis at a point $(\tau_{i},0)$, we have from~\eqref{hmep1}-\eqref{hmep3}, \eqref{c1}, \eqref{Casp4} and~\eqref{Casp5} that
\begin{align}
|\hat{m}_{i}(t(x),x)|\leq&|\hat{m}_{il}(\tau_{i})|+C\int_{0}^{x}|\hat{m}(t(y),y)|dy. \label{hapri1}
\end{align}
\indent For~\eqref{b13}, applying the same procedures as above, we further use~\eqref{Casp6} to obtain
\begin{align}
|\hat{n}_{i}(t(x),x)|\leq&|\hat{n}_{il}(\tau_{i})|+C(\int_{0}^{x}|\hat{n}(t(y),y)|dy+
\int_{0}^{x}|\hat{m}(t(y),y)|dy). \label{hapri2}
\end{align}
Taking the summation of~\eqref{hapri1} and~\eqref{hapri2} and the summation for $i=1,2$, applying the Gronwall's inequality, we have
\begin{align}\label{c6}
\max_{(t,x)\in E_{2}}|\hat{m}(t,x)|+|\hat{n}(t,x)|\leq C(\|\hat{m}_{l}\|_{C^{0}([0,+\infty))}+\|\hat{n}_{l}\|_{C^{0}([0,+\infty))}),\end{align}
where we have used the arbitrariness of $(t,x)\in E_{2}$.\\
\textbf{Region 3:} in the remaining region
$$
E_{3}=\{(t,x)|0\leq t\leq T_{1}, 0\leq x\leq L, x_{1}^{*}(t)\leq x \leq x_{2}^{*}(t)\}.
$$
\indent For any point $(t,x)\in E_{3}$, integrating the first equation in ~\eqref{b10} and~\eqref{b11} along the first characteristic curve that intersects the $x_{2}^{*}(t)$ at a point $(t_{1},x_{1})$, we get from~\eqref{mep1}-\eqref{mep3}, \eqref{nep1}, \eqref{nep2}, \eqref{c1}, \eqref{Casp4} and~\eqref{Casp5} that
\begin{align}
|m_{1}(t,x(t))|\leq&|m_{1}(t_{1},x_{1})|+C\int_{t_{1}}^{t}|m(\tau,x(\tau))|d\tau\notag\\
\leq&|m_{1}(t_{1},x_{1})|+C\int_{0}^{t}|m(\tau,x(\tau))|d\tau, \label{c7}\\
|n_{1}(t,x(t))|\leq&|n_{1}(t_{1},x_{1})|+C(\int_{0}^{t}|n(\tau,x(\tau))|d\tau
+\int_{0}^{t}|m(\tau,x(\tau))|d\tau). \label{c8}
\end{align}
\indent Similarly, for any point $(t,x)\in E_{3}$, integrating the second equation in~\eqref{b10} and~\eqref{b11} along the second characteristic curve that intersects $x_{1}^{*}(t)$ at a point $(t_{2},x_{2})$, we have
\begin{align}
|m_{2}(t,x(t))|\leq&|m_{2}(t_{2},x_{2})|+C\int_{0}^{t}|m(\tau,x(\tau))|d\tau,\label{c9}\\
|n_{2}(t,x(t))|\leq&|n_{2}(t_{2},x_{2})|+C(\int_{0}^{t}|n(\tau,x(\tau))|d\tau+\int_{0}^{t}|m(\tau,x(\tau))|d\tau). \label{c10}
\end{align}
\indent By applying the Gronwall's inequality, the combination of~\eqref{c7}-\eqref{c10} gives rise to
\begin{equation}\label{c11}
\begin{split}
\max_{(t,x)\in E_{3}}(|m(t,x)|+|n(t,x)|)\leq& C(\|m_{0}\|_{C^{0}([0,L])}+\|n_{0}\|_{C^{0}([0,L])}\\
&+\|\hat{m}_{l}\|_{C^{0}([0,+\infty))}+\|\hat{n}_{l}\|_{C^{0}([0,+\infty))}),
\end{split}
\end{equation}
where we have used~\eqref{Apri1} and~\eqref{c6} and the arbitrariness of $(t,x)\in E_{3}$.\\
\indent We notice from~\eqref{Apri1}, \eqref{c6}, \eqref{c11}, \eqref{b8} and~\eqref{hsdec1} that under the initial and boundary conditions~\eqref{a19}-\eqref{a20} for a sufficiently small $\varepsilon>0$ and the assumption~\eqref{Casp5}, we can check the validity of hypothesis ~\eqref{c1} for some constant $C>0$. Therefore, we obtain a uniform a-priori estimate for the global $C^{1}$ solution. The global existence of solution to the initial-boundary value problem~\eqref{b5} and~\eqref{1in1}-\eqref{1in2} can be checked by the standard continuity method, the details are omitted here.

\section{Periodic Solution}\label{s4}
\indent\indent In this section, we will prove global solution $V=(\rho(t,x),u(t,x))^{\top}$ is a time-periodic function with a period $P>0$.\\
\indent Using a Riemann invariant of system ~\eqref{a1}
\begin{align}\label{d1}
r=\frac{1}{2}(u-\frac{2}{\gamma-1}c),~~~s=\frac{1}{2}(u+\frac{2}{\gamma-1}c),
\end{align}
~\eqref{a1} can be converted into the following form
\begin{align}\label{d2}
\left\{
\begin{aligned}
r_{t}+\lambda_{1}(r,s)r_{x}=\frac{\beta(r+s)^{\alpha+1}}{2},\\
s_{t}+\lambda_{2}(r,s)s_{x}=\frac{\beta(r+s)^{\alpha+1}}{2},
\end{aligned}\right.
\end{align}
where
\begin{align*}
\lambda_{1}=u-c=\frac{\gamma+1}{2}r-\frac{\gamma-3}{2}s,\quad
\lambda_{2}=u+c=\frac{3-\gamma}{2}r+\frac{\gamma+1}{2}s.
\end{align*}
Correspondingly, the initial data and boundary conditions become
\begin{align}
&r(0,x)=r_{0}(x),~~s(0,x)=s_{0}(x),~~x\in[0,L]\label{d3},\\
&r(t,0)=r_{l}(t),~~~~s(t,0)=s_{l}(t),~~~\quad t\geq0\label{d4},
\end{align}
where $r_{l}(t),s_{l}(t)$ are time-periodic with the period $P>0$.\\
\indent For the convenience of later proof, we exchange $t$ and $x$, then problem ~\eqref{d2} and ~\eqref{d3}-~\eqref{d4} becomes the following Cauchy problem in the domain $E$
\begin{equation}\label{d5}
\left\{\begin{aligned}
&r_{x}+\frac{1}{\lambda_{1}}r_{t}=\frac{\beta(r+s)^{\alpha+1}}{2\lambda_{1}},\\
&s_{x}+\frac{1}{\lambda_{2}}s_{t}=\frac{\beta(r+s)^{\alpha+1}}{2\lambda_{2}},\\
&r(t,0)=r_{l}(t),\\
&s(t,0)=s_{l}(t).
\end{aligned}\right.
\end{equation}
Furthermore, setting
$$
W=(r-\tilde{r},s-\tilde{s})^{\top},\quad
\Lambda(t,x)=\left(
            \begin{array}{cc}
            \frac{1}{\lambda_{1}(r(t,x),s(t,x))} & 0\\
            0 & \frac{1}{\lambda_{2}(r(t,x),s(t,x))}\\
            \end{array}
\right),
$$
then ~\eqref{d5} can be rewritten as
\begin{align}\label{d6}
W_{x}+\Lambda(t,x) W_{t}=\frac{\beta}{2}\Lambda(t,x)\left(
                   \begin{aligned}
                   (r+s)^{\alpha+1}\\
                   (r+s)^{\alpha+1}\\
                   \end{aligned}
\right)
-\frac{\beta}{2}
\left(
            \begin{aligned}
            \frac{(\tilde{r}+\tilde{s})^{\alpha+1}}{\tilde{\lambda}_{1}}\\
            \frac{(\tilde{r}+\tilde{s})^{\alpha+1}}{\tilde{\lambda}_{2}}\\
            \end{aligned}
\right),
\end{align}
where
$$
\tilde{r}=\frac{1}{2}(\tilde{u}-\frac{2}{\gamma-1}\tilde{c}),\quad
\tilde{s}=\frac{1}{2}(\tilde{u}+\frac{2}{\gamma-1}\tilde{c}),
$$
\begin{align*}
&\tilde{\lambda}_{1}=\lambda_{1}(\tilde{r},\tilde{s})=\frac{\gamma+1}{2}\tilde{r}-\frac{\gamma-3}{2}\tilde{s},\\
&\tilde{\lambda}_{2}=\lambda_{2}(\tilde{r},\tilde{s})=\frac{3-\gamma}{2}\tilde{r}+\frac{\gamma+1}{2}\tilde{s}.
\end{align*}
\indent By
\begin{align*}
\|\rho-\tilde{\rho}\|_{C^{1}(E)}+\|u-\tilde{u}\|_{C^{1}(E)}<C\varepsilon,
\end{align*}
and~\eqref{d1}, we can get
\begin{align}\label{d7}
\|r(t,x)-\tilde{r}(x)\|_{C^{1}(E)}+\|s(t,x)-\tilde{s}(x)\|_{C^{1}(E)}<K_{1}\varepsilon
\end{align}
with $K_{1}>0$ a constant that depending only on $\tilde{\rho}, \tilde{u}, \gamma$ and $L$.\\
\indent Next we will show that the following conclusion holds
\begin{align}\label{d8}
r(t+P,x)=r(t,x),~~s(t+P,x)=s(t,x),\quad\forall t>T_{1},x\in[0,L],
\end{align}
where $T_{1}$ is defined by~\eqref{a21}.\\
\indent Letting
$$
U(t,x)=W(t+P,x)-W(t,x),
$$
then by ~\eqref{d6}, we can get
\begin{align}\label{d9}
\left\{\begin{aligned}
&U_{x}+\Lambda(t,x)U_{t}=G(t,x),\\
&U(t,0)=0,\quad t>0,
\end{aligned}\right.
\end{align}
where
\begin{align*}
G(t,x)=&\frac{\beta}{2}\Lambda(t+P,x)\left(
\begin{aligned}
(r(t+P,x)+s(t+P,x))^{\alpha+1}\\
(r(t+P,x)+s(t+P,x))^{\alpha+1}\\
\end{aligned}
\right)\\
&-\frac{\beta}{2}\Lambda(t,x)\left(
                            \begin{aligned}
                            (r(t,x)+s(t,x))^{\alpha+1}\\
                            (r(t,x)+s(t,x))^{\alpha+1}\\
                            \end{aligned}
\right)\\
&-[\Lambda(t+P,x)-\Lambda(t,x)]W_{t}(t+P,x).
\end{align*}
\indent Noting that $\lambda_{1}, \lambda_{2}$ are continuous functions of $(r,s)$, then by ~\eqref{d7}, we can get the following estimates
\begin{align}
&|W_{t}(t+p,x)|\leq K_{1}\varepsilon,\label{d10}\\
&|r(t+P,x)+s(t+P,x)|\leq K_{2},\label{d11}\\
&|\Lambda_{t}(r(t,x),s(t,x))|\leq K_{3}\varepsilon,\label{d12}\\
&|\Lambda(t+P,x)-\Lambda(t,x)|\leq K_{4}|U(t,x)|,\label{d13}\\
&|\Lambda(t,x)|\leq K_{5},\label{d14}
\end{align}
where constants $K_{2}, K_{3}, K_{4}, K_{5}$ depend only on $\tilde{\rho}, \tilde{u}, \gamma$ and $L$.\\
\indent It follows from ~\eqref{d10}-\eqref{d11}, \eqref{d13}-\eqref{d14} that
\begin{equation}\label{d15}
\begin{split}
|G(t,x)|\leq&\frac{|\beta|}{2}|\Lambda(t,x)|\left(
            \begin{aligned}
            (\alpha+1)|\eta|^{\alpha}|U(t,x)|\\
            (\alpha+1)|\eta|^{\alpha}|U(t,x)|\\
            \end{aligned}
\right)\\
&+\frac{|\beta|}{2}|\Lambda(t+P,x)-\Lambda(t,x)|\left(
            \begin{aligned}
            |r(t+P,x)+s(t+P,x)|^{\alpha+1}\\
            |r(t+P,x)+s(t+P,x)|^{\alpha+1}\\
            \end{aligned}
\right)\\
&+|\Lambda(t+P,x)-\Lambda(t,x)||W_{t}(t+P,x)|\\
\leq& K_{6}|U(t,x)|,
\end{split}
\end{equation}
where $\eta$ lies between $u(t,x)$ and $u(t+p,x)$, the definition of $K_{6}$ is the same as above.\\
\indent For a fixed point $(t_{0},x_{0})$ with $t_{0}>T_{1}, 0<x_{0}<L$, we can draw two  characteristic curves $\Gamma_{1}:t=t_{1}^{*}(x)$ and $\Gamma_{2}:t=t_{2}^{*}(x)$, namely,
\begin{align*}
\frac{dt_{1}^{*}}{dx}=\frac{1}{\lambda_{1}(r(t_{1}^{*},x),s(t_{1}^{*},x))},t_{1}^{*}(x_{0})=t_{0}
\end{align*}
and
\begin{align*}
\frac{dt_{2}^{*}}{dx}=\frac{1}{\lambda_{2}(r(t_{2}^{*},x),s(t_{2}^{*},x))},t_{2}^{*}(x_{0})=t_{0}
\end{align*}
for $0<x<x_{0}$. And we can easily see that $\Gamma_{1}$ lies below $\Gamma_{2}$.\\
Setting
\begin{align}\label{d16}
I(x)=\frac{1}{2}\int_{t_{1}^{*}(x)}^{t_{2}^{*}(x)}|U(t,x)|^{2}dt,
\end{align}
where $0\leq x< x_{0}$.\\
\indent By the definition of $T_{1}$ and $t_{0}>T_{1}$, we can get that $(t_{1}^{*}(0),t_{2}^{*}(0))\subset(0,+\infty)$, then by ~\eqref{d9}, we have $U(t,0)\equiv0$ in this interval.\\
\indent Therefore,
\begin{align*}
I(0)=0.
\end{align*}
\indent Taking derivative of $I(x)$ with respect to $x$, we get
\begin{align*}
I^{'}(x)=& \int_{t_{1}^{*}(x)}^{t_{2}^{*}(x)}{U(t,x)^{T}U_{x}(t,x)}dt+\frac{1}{2}|{U(t_{2}^{*}(x),x)}|^{2}
{\frac{1}{\lambda_{2}(r(t_{2}^{*}(x),x),s(t_{2}^{*}(x),x))}}\\
&-\frac{1}{2}|{U(t_{1}^{*}(x),x)}|^{2}{\frac{1}{\lambda_{1}(r(t_{1}^{*}(x),x),s(t_{1}^{*}(x),x))}}\\
\leq&-\int_{t_{1}^{*}(x)}^{t_{2}*(x)}U(t,x)^{T}\Lambda(t,x)U_{t}(t,x)dt
+\int_{t_{1}^{*}(x)}^{t_{2}^{*}(x)}U(t,x)^{T}G(t,x)dt\\
&+\frac{1}{2}U(t,x)^{T}\Lambda(t,x)U(t,x)|_{t=t_{1}^{*}(x)}^{t=t_{2}^{*}(x)}\\
=&-\frac{1}{2}\int_{t_{1}^{*}(x)}^{t_{2}^{*}(x)}(U(t,x)^{T}\Lambda(t,x)U(t,x))_{t}-U(t,x)^{T}
\Lambda_{t}(t,x)U(t,x)dt\\
&+\int_{t_{1}^{*}(x)}^{t_{2}^{*}(x)}U(t,x)^{T}G(t,x)dt+\frac{1}{2}U(t,x)^{T}\Lambda(t,x)U(t,x)
|_{t=t_{1}^{*}(x)}^{t=t_{2}^{*}(x)}\\
=&\frac{1}{2}\int_{t_{1}^{*}(x)}^{t_{2}^{*}(x)}U(t,x)^{T}\Lambda_{t}(t,x)U(t,x)dt
+\int_{t_{1}^{*}(x)}^{t_{2}^{*}(x)}U(t,x)^{T}G(t,x)dt\\
\leq&(K_{3}\varepsilon+2K_{6})I(x).
\end{align*}
In the last inequality we have used (\ref{d12}) and (\ref{d15}).\\
\indent Hence, by Gronwall's inequality, we can get that $I(x)\equiv0$.
Furthermore, by continuity of $I(x)$, we have $I(x_{0})=0$, then $U(t_{0},x_{0})=0$.\\
\indent Since $(t_{0},x_{0})$ is arbitrary, so we have
$$
U(t,x)\equiv0,\quad \forall t>T_{1},x\in[0,L],
$$
 that is,  we complete the proof of ~\eqref{d8}. Then, using~\eqref{d1} and $c=\sqrt{\gamma}\rho^{\frac{\gamma-1}{2}}$, we can get $(\rho,u)^{\top}$ is also a periodic function with a period $P>0$.

\end{sloppypar}
\end{document}